\documentclass[10pt]{article}

\usepackage{times}
\usepackage{graphicx} 
\usepackage{subfigure} 

\usepackage{natbib}

\usepackage{amsmath,amsfonts,amsthm} 
\usepackage{amsmath, amssymb, color}

\usepackage{algorithm}
\usepackage{algorithmic}

\usepackage{hyperref}

\usepackage[accepted]{icml2017}

\newtheorem{assumption}{Assumption}

\newtheorem{theorem}{Theorem}
\newtheorem{lemma}[theorem]{Lemma}

\newcommand{\bb}{b}

\newcommand{\real}[1]{\mbox{$\mathbb{R}^{#1}$}}

\newcommand{\Rmnum}[1]{\uppercase\expandafter{\romannumeral #1}} 
\usepackage[raggedright]{titlesec}
\titleformat{\chapter}{\centering\Huge\bfseries}{Chapter \Rmnum{\thechapter} }{1em}{} 

\usepackage{mathrsfs} 
\usepackage{enumerate} 
\graphicspath{{figures/}}  %

\usepackage{algorithm}
\usepackage{algorithmic}
\usepackage{multirow}
\usepackage{amsmath}
\usepackage{stmaryrd}
\usepackage{float}
\usepackage{graphicx}
\usepackage{subfigure}
\usepackage{xcolor}
\usepackage{color}

\begin{document}

\twocolumn[
\icmltitle{Strong NP-Hardness for Sparse Optimization with Concave Penalty Functions}



\icmlsetsymbol{equal}{*}

\begin{icmlauthorlist}
\icmlauthor{Yichen Chen}{pu}
\icmlauthor{Dongdong Ge}{sh}
\icmlauthor{Mengdi Wang}{pu}
\icmlauthor{Zizhuo Wang}{um}
\icmlauthor{Yinyu Ye}{stanford}
\icmlauthor{Hao Yin}{stanford}
\end{icmlauthorlist}

\icmlaffiliation{pu}{Princeton University, NJ, USA}
\icmlaffiliation{sh}{Shanghai University of Finance and Economics, Shanghai, China}
\icmlaffiliation{um}{University of Minnesota, MN, USA}
\icmlaffiliation{stanford}{Stanford University, CA, USA}

\icmlcorrespondingauthor{Mengdi Wang}{mengdiw@princeton.edu}

\icmlkeywords{boring formatting information, machine learning, ICML}
\vskip 0.3in
]

\printAffiliationsAndNotice{}

%

\begin{abstract}
Consider the regularized sparse minimization problem, which involves empirical sums of loss functions for $n$ data points (each of dimension $d$) and a nonconvex sparsity penalty.  We prove that finding an $\mathcal{O}(n^{c_1}d^{c_2})$-optimal solution to the regularized sparse optimization problem is strongly NP-hard for any $c_1, c_2\in [0,1)$ such that $c_1+c_2<1$. The result applies to a broad class of loss functions and sparse penalty functions.  It suggests that one cannot even approximately solve the sparse optimization problem in polynomial time, unless P $=$ NP.
\end{abstract}


\section{Introduction}

We study the sparse minimization problem, where the objective is the sum of empirical losses over input data and a sparse penalty function. Such problems commonly arise from empirical risk minimization and variable selection. 
The role of the penalty function is to induce sparsity in the optimal solution, i.e., to minimize the empirical loss using  as few nonzero coefficients as possible. 

{\bf Problem 1} Given the loss function $\ell:\mathbb{R}\times \mathbb{R}\mapsto \mathbb{R}^+$, penalty function $p:\mathbb{R}\mapsto \mathbb{R}^+$, and regularization parameter $\lambda >0$, consider the problem 
\begin{equation*}
  \min_{x\in \mathbb{R}^d} \sum_{i=1}^n \ell\left(a_i^T x,b_i\right)+\lambda \sum_{j=1}^d p\left(|x_j|\right) ,
 \end{equation*}
where $ A=(a_1,\ldots,a_n)^T$ $\in \mathbb{R}^{n\times d}$, $b=(b_1,\ldots,b_n)^T$ $\in \mathbb{R}^n$ are input data.



%

%

We are interested in the computational complexity of Problem 1 under general conditions of the loss function $\ell$ and the sparse penalty $p$. 
In particular, we focus on the case where $\ell$ is a convex loss function and $p$ is a concave penalty with a unique minimizer  at $0$. Optimization problems with convex $\ell$ and concave $p$ are common in sparse regression, compressive sensing, and sparse approximation. A list of applicable examples of $\ell$ and $p$ is given in  Section 3.

For certain special cases of Problem 1, it has been shown that finding an {\it exact solution} is strongly NP-hard \citep{huo2010complexity, ChenX_GeD_2014}. 
However, these results have not excluded the possibility of the existence of polynomial-time algorithms with small approximation error. 
\cite{ChenY_2016} established the hardness of approximately solving Problem 1 when $p$ is the $L_0$ norm.

%

In this paper, 
we prove that it is strongly NP-hard to approximately solve Problem 1 within certain optimality error. 
More precisely, we show that there exists a lower bound on the suboptimality error of any polynomial-time deterministic algorithm.  
Our results apply to a variety of optimization problems in estimation and machine learning. Examples include sparse classification, sparse logistic regression, and many more. 
The strong NP-hardness of approximation is one of the strongest forms of complexity result for continuous optimization. To our best knowledge, this paper gives the first and strongest set of hardness results for Problem 1 under very general assumptions regarding the loss and penalty functions.

Our main contributions are three-fold.
\begin{enumerate}
\item We prove the strong NP-hardness for Problem 1 with general loss functions. This is the first results that apply to the broad class of problems including but not limited to: least squares regression, linear model with Laplacian noise, robust regression, Poisson regression, logistic regression, inverse Gaussian models, etc.

\item We present a general condition on the sparse penalty function $p$ such that Problem 1 is strongly NP-hard. The condition is a slight weaker version of strict concavity. It is satisfied by typical penalty functions such as the $L_q$ norm ($q\in[0,1)$), clipped $L_1$ norm, SCAD, etc.  To the best of our knowledge, this is the most general condition on the penalty function in the literature.

\item We prove that finding an $\mathcal{O}\left(\lambda n^{c_1} d^{c_2}\right)$-optimal solution to Problem 1 is strongly NP-hard, for any $c_1, c_2\in [0,1)$ such that $c_1+c_2<1$. Here the $\mathcal{O}(\cdot)$ hides parameters that depend on the penalty function $p$, which is to be specified later. 
It illustrates a gap between the optimization error achieved by any tractable algorithm and the desired statistical precision. 
Our proof provides a first unified analysis that deals with a broad class of problems taking the form of Problem 1. 


\end{enumerate}

Section 2 summarizes related literatures from optimization, machine learning and statistics. Section 3 presents the key assumptions and illustrates examples of loss and penalty functions that satisfy the assumptions. Section 4 gives the main results. Section 5 discusses the implications of our hardness results. Section 6 provides a proof of the main results in a simplified setting.  
The full proofs are deferred to the appendix.

\section{Background and Related Works}
\label{Sec:Literature}

Sparse optimization is a powerful machine learning tool for extracting useful information for massive data. In Problem 1, the sparse penalty serves to select the most relevant variables from a large number of variables, in order to avoid overfitting. In recent years, nonconvex choices of $p$ have received much attention; see \cite{frank1993statistical,fan2001variable,chartrand2007exact, candes2008enhancing,  fan2010selective, xue2012nonconcave, loh2013regularized, wang2014optimal, fan2015tac}.

%
%

Within the optimization and mathematical programming community, the complexity of Problem 1 has been considered in a number of special cases.  \cite{huo2010complexity} first proved the hardness result for a relaxed family of penalty functions with $L_2$ loss. They show that for the penalties in $L_0$, hard-thresholded \citep{antoniadis2001regularization} and SCAD \citep{fan2001variable}, the above optimization problem is NP-hard.  \cite{ChenX_GeD_2014} showed that the $L_2$-$L_p$ minimization is strongly NP-hard when $p\in(0,1)$. At the same time, \citep{bian2014optimality} proved the strongly NP-hardness for another class of penalty functions. 
The preceding existing analyses mainly focused on finding an exact global optimum to Problem 1. For this purpose, they implicitly assumed that all the input and parameters involved in the reduction are rational numbers with a finite numerical representation, otherwise finding a global optimum to a continuous problem would be always intractable. 
A recent technical report \cite{ChenY_2016} proves the hardness of obtaining an $\epsilon$-optimal solution when $p$ is the $L_0$ norm. 

Within the theoretical computer science community, there have been several early works on the complexity of sparse recovery, beginning with \cite{arora1993hardness}.  \cite{amaldi1998approximability} proved that 
the problem $\min \{\|x\|_0 \mid  Ax=b\}$ is not approximable within a factor $2^{\log^{1-\epsilon}d}$ for any $\epsilon>0$. \cite{natarajan1995sparse}  showed that, given $\epsilon>0, A$ and $b$, the problem $\min\{ \|x\|_0 \mid \|Ax-b\|_2\leq \epsilon\}$ is NP-hard. \cite{davis1997adaptive} proved a similar result that for some given $\epsilon>0$ and $M>0$, it is NP-complete to find a solution $x$ such that $\|x\|_0\leq M$ and $\|Ax-b\|\leq \epsilon$.
More recently, \cite{foster2015variable} studied sparse recovery and sparse linear regression with subgaussian noises. Assuming that the true solution is $K$-sparse, it showed that no polynomial-time (randomized) algorithm can find a $K\cdot 2^{\log^{1-\delta} d}$-sparse solution $x$ with $ \|Ax-b\|_2^2\leq d^{C_1} n^{1-C_2}$ {\it with high probability}, where $\delta, C_1, C_2$ are arbitrary positive scalars. 
Another work \cite{zhang2014lower} showed that under the Gaussian linear model, there exists a gap between the mean square loss that can be achieved by polynomial-time algorithms and the statistically optimal mean squared error. These two works focus on estimation of linear models and impose distributional assumptions regarding the input data. These results on estimation are different in nature with our results on optimization.

In contrast, we focus on the optimization problem itself. Our results apply to a variety of loss functions and penalty functions, not limited to linear regression. Moreover, we do not make any distributional assumption regarding the input data.

There remain several open questions. First, existing results mainly considered least square problems or $L_q$ minimization problems. Second, existing results focused mainly on the $L_0$ penalty function. The complexity of Problem 1 with general loss function and penalty function is yet to be established. Things get complicated when $p$ is a continuous function instead of the discrete $L_0$ norm function. The complexity for finding an $\epsilon$-optimal solution with general $\ell$ and $p$ is not fully understood. We will address these questions in this paper.

\section{Assumptions}

In this section, we state the two critical assumptions that lead to the strong NP-hardness results: one for the penalty function $p$, the other one for the loss function $\ell$. We argue that these assumptions are essential and very general. They apply to a broad class of loss functions and penalty functions that are commonly used.

\subsection{Assumption About Sparse Penalty}

Throughout this paper, we make the following assumption regarding the sparse penalty function $p(\cdot)$.

\begin{assumption}
\label{asm:pen}
The penalty function
$p(\cdot)$ satisfies the following conditions:
\begin{enumerate} [(i)]
\item (Monotonicity) $p(\cdot)$ is non-decreasing on $[0,+\infty)$.
\item {(Concavity)} There exists $\tau >0$ such that $p(\cdot)$ is concave but not linear on $[0, \tau]$.
\end{enumerate}
\end{assumption}

    In words, condition (ii) means that the concave penalty $p(\cdot)$ is nonlinear. Assumption \ref{asm:pen} is the most general condition on penalty functions in the existing literature of sparse optimization. Below we present a few such examples. 

\begin{enumerate}

\item In variable selection problems, the $L_0$ penalization
$p(t)=I_{\{t\neq 0\}}$ arises naturally as a penalty for the number of factors selected.

\item A natural generalization of the $L_0$ penalization is the $L_p$ penalization $p(t)=t^p$ where $(0< p <1)$. The corresponding minimization problem is called the bridge regression problem \cite{frank1993statistical}.

\item To obtain a hard-thresholding estimator, \citet{antoniadis2001regularization} use the penalty functions $p_\gamma (t)=\gamma^2-((\gamma-t)^+)^2$ with $\gamma >0$, where $(x)^+ := \max\{x,0 \}$ denotes the positive part of $x$.

\item Any penalty function that belongs to the folded concave penalty family \cite{FanJ_XueL_2014} satisfies the conditions in Theorem \ref{Thm_Main}. Examples include the SCAD \cite{fan2001variable} and the MCP \cite{ZhangCH_2010}, whose derivatives on $(0, +\infty)$ are
$p^\prime_\gamma(t)=\gamma
I_{\{t\leq\gamma\}}+\frac{(a\gamma-t)^+}{a-1} I_{\{t>\gamma\}}
\quad \text{and} \quad p^\prime_\gamma(t)=(\gamma-\frac{t}{b})^+, $ respectively, where $\gamma>0$, $a>2$ and $b>1$.

\item The conditions in Theorem \ref{Thm_Main} are also satisfied by the clipped $L_1$ penalty function \cite{antoniadis2001regularization,ZhangT_2010} $p_\gamma(t)=\gamma \cdot \min(t, \gamma)$ with $\gamma>0$. This is a special case of the piecewise linear penalty function:
\begin{equation*} \label{Eq:PiecewiseLinear}
p(t)=\left\{
\begin{array}{ll}
k_1t&\mbox{if $0 \leq t \leq a$}\\
k_2t+(k_1-k_2)a & \mbox{if $t> a$}
\end{array}
\right.
\end{equation*}
where $0\leq k_2<k_1$ and $a>0$.

\item Another family of penalty functions which bridges the $L_0$ and $L_1$ penalties are the fraction penalty functions $\displaystyle p_\gamma (t)=\frac{(\gamma+1)t}{\gamma+t}$ with $\gamma>0$ \cite{LvJ_FanY_2009}.

\item The family of log-penalty functions:
$$
p_\gamma(t)=\frac{1}{\log (1+\gamma)}\log(1+\gamma t)
$$
with $\gamma >0$, also bridges the $L_0$ and $L_1$ penalties \cite{candes2008enhancing}.%

\end{enumerate}

\subsection{Assumption About Loss Function}
\label{sec:app}

We state our assumption about the loss function $\ell$. 

\begin{assumption}
 \label{asm:loss}
Let $M$ be an arbitrary constant. 
For any interval $[\tau_1,\tau_2]$ where $0<\tau_1<\tau_2<M$, there exists $k\in \mathbb{Z}^+$ and $b \in \mathbb{Q}^{k}$ such that $h(y)=\sum_{i=1}^{k} \ell(y,b_i)$ has the following properties: 
\begin{enumerate}[(i)]
\item $h(y)$ is convex and Lipschitz continuous on $[\tau_1,\tau_2]$.
\item $h(y)$ has a unique minimizer  $y^*$ in $(\tau_1,\tau_2)$.
\item \label{asm:iii} There exists $N\in \mathbb{Z^+}, \bar{\delta}\in \mathbb{Q^+}$ and $C\in \mathbb{Q^+}$ such that when $\delta\in (0,\bar{\delta})$, we have 
$$ \frac{h(y^*\pm\delta)-h(y^*)}{\delta^N}\geq C.$$
\item $h(y^*)$, $\{b_i\}^k_{i=1}$ can be represented in $\mathcal{O}( \log\frac1{\tau_2-\tau_1})$ bits. 
\end{enumerate}
\end{assumption}

Assumption \ref{asm:loss} is a critical, but very general,  assumption regarding the loss function $\ell(y,b)$. 
Condition (i) requires convexity and Lipschitz continuity within a neighborhood.
Conditions (ii), (iii) essentially require that, given an interval $[\tau_1,\tau_2]$, one can artificially pick $b_1,\ldots,b_k$ to construct a function $h(y)=\sum_{i=1}^{k} \ell(y,b_i)$ such that $h$ has its unique minimizer in $[\tau_1,\tau_2]$ and has enough curvature near the minimizer. This property ensures that a bound on the minimal value of $h(y)$ can be translated to a meaningful bound on the minimizer $y^*$.
The conditions (i), (ii), (iii) are typical properties that a loss function usually satisfies. 
Condition (iv) is a technical condition that is used to avoid dealing with infinitely-long irrational numbers. It can be easily verified for almost all common loss functions. 
 
We will show that Assumptions 2 is satisfied by a variety of loss functions. An (incomplete) list is given below.

\begin{enumerate}
\item In the least squares regression, the loss function has the form  
\begin{equation*}
\begin{aligned}
\sum_{i=1}^n \left(a_i^T x-b_i\right)^2.
\end{aligned}
\end{equation*}
Using our notation, the corresponding loss function is $\ell(y,b)=(y-b)^2$. For all $\tau_1,\tau_2$, we choose an arbitrary $b^\prime\in [\tau_1,\tau_2]$. We can verify that $h(y)=\ell(y,b^\prime)$ satisfies all the conditions in Assumption \ref{asm:loss}. 

\item In the linear model with Laplacian noise,  the negative log-likelihood function is 
\begin{equation*}
  \sum_{i=1}^n \left|a_i^T x-b_i\right|.
\end{equation*}
So the loss function is $\ell(y,b)=|y-b|$. As in the case of least squares regression, the loss function satisfy Assumption \ref{asm:loss}. Similar argument also holds when we consider the $L_q$ loss $|\cdot|^q$ with $q\geq 1$.

\item In robust regression, we consider the Huber loss \cite{huber1964robust} which is a mixture of $L_1$ and $L_2$ norms. 
The loss function takes the form
\begin{equation*} 
 L_\delta(y,b )=\left\{\begin{array}{ll}
                     \frac{1}{2}|y-b|^2 &\text{ for $|y -b|\leq \delta$},\\
                     \delta(|y -b|-\frac{1}{2}\delta) &\text{ otherwise.}\\
                    \end{array}\right.
\end{equation*}
for some $\delta>0$ where $y = a^Tx$. We then verify that Assumption \ref{asm:loss} is satisfied. For any interval $[\tau_1, \tau_2]$, we pick an arbitrary $b \in [\tau_1,\tau_2]$ and let $h(y)=\ell(y,b)$. We can see that $h(y)$ satisfies all the conditions in Assumption \ref{asm:loss}. 

\item In Poisson regression \cite{cameron2013regression}, the negative log-likelihood minimization is
\begin{equation*}
    \underset{x\in \mathbb{R}^d}{\text{min}} {-\log L(x;A,b)} =\underset{x\in \mathbb{R}^d}{\text{min}}\sum_{i=1}^n (\exp(a_i^Tx)-b_i \cdot a_i^Tx ).
\end{equation*}
We now show that  $\ell(y,b)=e^y-b\cdot y$ satisfies Assumption \ref{asm:loss}. For any interval $[\tau_1,\tau_2]$, we choose $q$ and $r$ such that $q/r\in[e^{\tau_1},e^{\tau_2}]$. Note that $e^{\tau_2}-e^{\tau_1}=e^{\tau_1+\tau_2-\tau_1}-e^{\tau_1}\geq \tau_2-\tau_1$. Also, $e^{\tau_2}$ is bounded by $e^{M}$.  Thus, $q,r$ can be chosen to be polynomial in $\lceil1/(\tau_2-\tau_1)\rceil$ by letting $r= \lceil1/(\tau_2-\tau_1)\rceil$ and $q$ be some number less than $r\cdot e^{M}$. Then, we choose $k=r$ and $b\in \mathbb{Z}^k$ such that $h(y)=\sum_{i=1}^k \ell(y,b_i)=r \cdot e^y-q\cdot y$. Let us verify Assumption \ref{asm:loss}. (i), (iv) are straightforward by our construction. For (ii), note that $h(y)$ take its minimum at $\ln (q/r)$ which is inside $[\tau_1,\tau_2]$ by our construction. To verify (iii), consider the second order Taylor expansion of $h(y)$ at $\ln (q/r)$,
\begin{align*}
    h(y+\delta)-h(y)= \frac{r\cdot e^y}{2} \cdot \delta^2 +o(\delta^2)\geq\frac{\delta^2}{2}+o(\delta^2),
\end{align*}
We can see that (iii) is satisfied. Therefore, Assumption \ref{asm:loss} is satisfied. 

\item In logistic regression, 
the negative log-likelihood function minimization is 
\begin{align*}
\underset{x\in \mathbb{R}^d}{\text{min}} \sum_{i=1}^n \log (1+\exp(a_i^T x))-\sum_{i=1}^n b_i\cdot a_i^T x.
\end{align*}
We claim that the loss function $\ell(y,b)=\log(1+\exp(y))-b\cdot y$ satisfies Assumption \ref{asm:loss}. By a similar argument as the one in Poisson regression, we can verify that $h(y)=\sum_{i=1}^r \ell (y,b_i)=r\log(1+\exp(y))-q y$ where $q/r\in [\frac{e^{\tau_1}}{1+e^{\tau_1}},\frac{e^{\tau_2}}{1+e^{\tau_2}}]$ and $q,r$ are polynomial in $\lceil1/(\tau_2-\tau_1)\rceil$ satisfies all the conditions in Assumption \ref{asm:loss}. For (ii), observe that $\ell(y,b)$ take its minimum at $y=\ln \frac{q/r}{1-q/r}$.  To verify (iii), we consider the second order Taylor expansion at $y=\ln\frac{q/r}{1-q/r}$, which is 
\begin{equation*}
    h(y+\delta)-h(y)=\frac{q}{2(1+e^y)}\delta^2+o(\delta^2)
\end{equation*}
where $y\in[\tau_1,\tau_2]$. Note that $e^y$ is bounded by $e^M$, which can be computed beforehand. As a result, (iii) holds as well. 

\item In the mean estimation of inverse Gaussian models \cite{mccullagh1984generalized}, the negative log-likelihood function minimization is
\begin{equation*}
    \underset{x\in \mathbb{R}^d}{\text{min}}\sum_{i=1}^n \frac{(b_i\cdot \sqrt{a_i^T x}-1)^2}{b_i}.
\end{equation*}
Now we show that the loss function $\ell(y,b)=\frac{(b\cdot \sqrt{y}-1)^2}{b}$ satisfies Assumption \ref{asm:loss}. By setting the derivative to be zero with regard to $y$, we can see that $y$ take its minimum at $y=1/b^2$. Thus for any $[\tau_1,\tau_2]$, we choose $b^\prime=q/r\in[1/\sqrt{\tau_2},1/\sqrt{\tau_1} ]$. We can see that $h(y)=\ell(y,b^\prime)$ satisfies all the conditions in Assumption \ref{asm:loss}. 
\item In the estimation of generalized linear model under the exponential distribution \cite{mccullagh1984generalized}, the negative log-likelihood function minimization is
\begin{equation*}
   \underset{x\in \mathbb{R}^d}{\text{min}}-\log L(x;A,b)= \underset{x\in \mathbb{R}^d}{\text{min}}\frac{b_i}{a_i^T x}+\log (a_i^T x).
\end{equation*}
By setting the derivative to 0 with regard to $y$, we can see that $\ell(y,b)=\frac{b}{y}+\log y$ has a unique minimizer at $y=b$. Thus by choosing $b^\prime\in [\tau_1,\tau_2]$ appropriately, we can readily show that $h(y)=\ell(y,b^\prime)$ satisfies all the conditions in Assumption \ref{asm:loss}. 
\end{enumerate}

To sum up, the combination of {\it any} loss function given in Section 3.1 and {\it any} penalty function given in Section 3.2 results in a strongly NP-hard optimization problem. 

\section{Main Results}

In this section, we state our main results on the strong NP-hardness of Problem 1. We warm up with a preliminary result for a special case of Problem 1.  

\begin{theorem}[A Preliminary Result]\label{Thm_Main}
Let Assumption \ref{asm:pen} hold, and let $p(\cdot)$ be twice continuously differentiable in $(0,\infty)$.
Then the minimization problem 
\begin{equation}
    \label{Eq:Problem_Main}
\min _{x \in \mathbb{R}^n}  ~ \|Ax-b\|_q^q + \lambda  \sum _{j=1}^d p(|x_j|),
\end{equation}
 is strongly NP-hard.
\end{theorem}
The result shows that many of the penalized least squares problems, e.g., \citep{fan2010selective}, while enjoying small estimation errors, are hard to compute. It suggests that there does not exist a fully polynomial-time approximation scheme for Problem 1. It has not answered the question: whether one can approximately solve Problem 1 within certain constant error.

Now we show that it is not even possible to efficiently approximate the global optimal solution of Problem 1, unless $P=NP$.  
Given an optimization problem $\min_{x\in X} f(x)$, we say that a solution $\bar x$ is $\epsilon$-optimal if $\bar x\in X$ and 
$f(\bar x)\leq \inf_{x\in X} f(x) +\epsilon.$

\begin{theorem}[Strong NP-Hardness of Problem 1]
    \label{thm:thm1}
Let Assumptions \ref{asm:pen} and \ref{asm:loss} hold, and let $c_1,c_2\in[0,1)$ be arbitrary such that $c_1+c_2< 1$. Then it is strongly NP-hard to find a $\lambda \cdot \kappa\cdot n^{c_1}d^{c_2}$-optimal solution of Problem 1, where $d$ is the dimension of variable space and $\kappa =\min_{t\in[\tau/2,\tau]} \{\frac{2p(t/2)-p(t)}{t}\}$.
\end{theorem}

The non-approximable error in Theorem \ref{thm:thm1} involves the constant 
$\kappa$ which is determined by the sparse penalty function $p$. 
In the case where $p$ is the $L_0$ norm function, we can take $\kappa  = 1$.  In the case of piecewise linear $L_1$ penalty, we have $\kappa = (k_1-k_2)/4 $. In the case of SCAD penalty, we have $\kappa = \Theta(\gamma^2).$

According to Theorem \ref{thm:thm1}, the non-approximable error $\lambda \cdot \kappa\cdot n^{c_1}d^{c_2}$ is determined by three factors: (i) properties of the regularization penalty $\lambda\cdot \kappa$; (ii) data size $n$; and (iii) dimension or number of variables $d$. This result illustrates a fundamental gap that can not be closed by any polynomial-time deterministic algorithm. This gap scales up when either the data size or the number of variables increases. In Section 5.1, we will see that this gap is substantially larger than the desired estimation precision in a special case of sparse linear regression.

\vspace{5pt}

Theorems \ref{Thm_Main} and \ref{thm:thm1} validate the long-lasting belief that optimization involving nonconvex penalty is hard. More importantly, Theorem \ref{thm:thm1} provide lower bounds for the optimization error that can be achieved by any polynomial-time algorithm.
This is one of the strongest forms of hardness result for continuous optimization.

\section{An Application and Remarks}

In this section, we analyze the strong NP-hardness results in the special case of linear regression with SCAD penalty (Problem 1). We give a few remarks on the implication of our hardness results.

\subsection{Hardness of Regression with SCAD Penalty}

Let us try to understand how significant is the non-approximable error of Problem 1. We consider the special case of linear models with SCAD penalty. 
Let the input data $(A,b)$ be generated by the linear model $A\bar x +\varepsilon =b$, where $\bar x$ is the unknown {\it true} sparse coefficients and $\varepsilon$ is a zero-mean multivariate subgaussian noise. Given the data size $n$ and variable dimension $d$, we follow \cite{fan2001variable} and obtain a special case of Problem 1, given by
\begin{equation}
    \min_x\frac{1}{2}\|Ax-b\|^2_2 + n \sum_{j=1}^d p_\gamma (|x_j|),
    \label{equ:scad}
\end{equation}
where $\gamma=\sqrt{\log d / n}$.
\cite{fan2001variable} showed that the optimal solution $x^*$ of problem \eqref{equ:scad} has  a small statistical error, i.e., $\|\bar x-x^*\|_2^2= \mathcal{O}\left( n^{-1/2} + a_n\right),$
where $a_n=\max \{p^\prime_\lambda(|x^*_{j}|):x^*_j\neq 0\}$. \cite{fan2015tac} further showed that we only need to find a $\sqrt{n \log d}$-optimal solution to \eqref{equ:scad} to achieve such a small estimation error.

However, Theorem 2 tells us that it is not possible to compute an $\epsilon_{d,n}$-optimal solution for problem \eqref{equ:scad} in polynomial time, where $\epsilon_{d,n} = \lambda \kappa n^{1/2}d^{1/3}$ (by letting $c_1=1/2,c_2=1/3$). In the special case of problem \eqref{equ:scad}, we can verify that $\lambda =n$ and $\kappa = \Omega(\gamma^2) =\Omega({\log d / n})$. As a result, we see that 
$$\epsilon_{d,n} = \Omega(n^{1/2} d^{1/3}) \gg \sqrt{n \log d},$$
for high values of the dimension $d$. 
According to Theorem 2, it is strongly NP-hard to approximately solve problem $\eqref{equ:scad}$ within the required statistical precision $\sqrt{n \log d}$.
This result illustrates a sharp contrast between statistical properties of sparse estimation and the worst-case computational complexity.

\subsection{Remarks on the NP-Hardness Results}

As illustrated by the preceding analysis, the non-approximibility of Problem 1 suggests that computing the sparse estimator is hard. The results suggest a fundamental conflict between computation efficiency and estimation accuracy in sparse data analysis. 
Although the results seem negative, they should not discourage researchers from studying computational perspectives of sparse optimization. We make the following remarks:
\begin{enumerate}
\item Theorems 1, 2 are {\it worst-case} complexity results. They suggest that one cannot find a tractable solution to the sparse optimization problems, without making any additional assumption to rule out the worst-case instances.
\item Our results do not exclude the possibility that, under more stringent modeling and distributional assumptions, the problem would be tractable with high probability or on average. 
\end{enumerate}
In short, the sparse optimization Problem 1 is fundamentally hard from a purely computational perspective.
This paper together with the prior related works provide a complete answer to the computational complexity of sparse optimization.

%
%
%
%
%

\section{Proof of Theorem \ref{Thm_Main}}\label{sec_prf}
In this section, we prove Theorem \ref{Thm_Main}. The proof of Theorems \ref{thm:thm1} 
is deferred to the appendix which is based on the idea of the proof in this section. 
We construct a polynomial-time reduction from the  {\it 3-partition problem} \cite{garey1978strong} to the sparse optimization problem. Given a set $S$ of $3m$ integers $s_1,...s_{3m}$, the three partition problem is to determine whether $S$ can be partitioned into $m$ triplets such that the sum of the numbers in each subset is equal. This problem is known to be strongly NP-hard \cite{garey1978strong}. 
The main proof idea bears a similar spirit as the works by \citet{huo2010complexity}, \citet{ChenX_GeD_2014} and  \citet{ChenY_2016}. 
The proofs of all the lemmas can be found in the appendix.

We first illustrate several properties of the penalty function if it
satisfies the conditions in Theorem \ref{Thm_Main}.

\begin{lemma}\label{Lem_PenaltyProperty}
If $p(t)$ satisfies the conditions in Theorem \ref{Thm_Main}, then 
for any $l \geq 2$, and any $t_1, t_2, \ldots, t_l \in \real{}$
, we have $p(|t_1|) + \cdots + p(|t_l|) \ge \min\{ p(|t_1 +\cdots + t_l|), p(\tau)\}$.
\end{lemma}

\begin{lemma}\label{Lem_PenaltyProperty2}
If $p(t)$ satisfies the conditions in Theorem \ref{Thm_Main}, then there exists
$\tau_0 \in (0, \tau)$ such that $p(\cdot)$ is concave but not linear on $[0, \tau_0]$ and is twice continuously differentiable on $[\tau_0, \tau]$. Furthermore, for any $\tilde t \in (\tau_0, \tau)$,
let $\bar \delta = \min\{ \tau_0 / 3, \tilde t - \tau_0, \tau - \tilde t\}$. Then for any $\delta \in (0, \bar \delta)$
$l \geq 2$, and any
$t_1, t_2, \ldots, t_l$ such that $t_1 +\cdots + t_l = \tilde t$, we have
$$
p(|t_1|) + \cdots + p(|t_l|) < p(\tilde t) + C_1 \delta
$$
only if $|t_i - \tilde{t}| < \delta$ for some $i$ while $|t_j| < \delta$ for all $j \neq i$, where $C_1 = \frac{p(\tau_0 / 3) + p(2 \tau_0 / 3) - p(\tau_0)}{\tau_0 / 3} > 0$.
\end{lemma}

In our proof of Theorem \ref{Thm_Main}, we will consider the
following function
\begin{equation*}\label{Eq:gDef}
g_{\theta, \mu}(t) := p(|t|)+ \theta \cdot |t|^q + \mu \cdot |t - \hat \tau|^q
\end{equation*}
with $\theta, \mu >0$, where $\hat \tau$ is an arbitrary fixed rational number in $(\tau_0, \tau)$. We have the following lemma about $g_{\theta, \mu}(t)$.

\begin{lemma}\label{Lem_MinUnique}
    If $p(t)$ satisfies the conditions in Theorem \ref{Thm_Main}, $q > 1$, and $\tau_0$ satisfies the properties in Lemma \ref{Lem_PenaltyProperty2}, then there exist $\underline \theta > 0$ and $\underline \mu > 0$ such that for any $\theta \geq \underline \theta$ and $\mu \geq \underline \mu \cdot \theta$,
the following properties are satisfied:
\begin{enumerate}  
\item $g^{\prime \prime}_{\theta, \mu}(t) \geq 1$ for any $t \in [\tau_0, \tau]$;
\item $g_{\theta, \mu}(t)$ has a unique global minimizer $t^*(\theta, \mu) \in (\tau_0, \tau)$;
\item Let $\bar \delta = \min\{t^*(\theta, \mu) - \tau_0, \tau - t^*(\theta, \mu), 1\}$, then for any $\delta \in (0, \bar \delta)$, we have $g_{\theta, \mu}(t) < h(\theta, \mu ) + \delta^2$ only if $|t - t^*(\theta, \mu)| < \delta$, where $h(\theta, \mu )$ is the minimal value of $g_{\theta, \mu}(t)$.
\end{enumerate}
\end{lemma}

\begin{lemma}\label{Lem_MinUnique_1}
If $p(t)$ satisfies the conditions in Theorem \ref{Thm_Main}, $q = 1$, and $\tau_0$ satisfies the properties in Lemma
\ref{Lem_PenaltyProperty2}, then there exist $\hat {\underline \mu} > 0$
 such that for any $\mu \geq \hat {\underline \mu}$, the following
properties are satisfied:
\begin{enumerate}  
\item $g_{0, \mu} ^\prime (t) < -1$ for any $t \in [\tau_0, \hat \tau)$ and $ g_{0, \mu} ^\prime (t) > 1$ for any $t \in (\hat \tau, \tau]$;
\item $ g_{0, \mu}(t)$ has a unique global minimizer $t^*(0, \mu)  = \hat \tau \in (\tau_0, \tau)$;
\item Let $\bar \delta = \min\{\hat \tau - \tau_0, \tau - \hat \tau, 1\}$, then for any $\delta \in (0, \bar \delta)$, we have $g_{0, \mu}(t) < h(0, \mu ) + \delta^2$ only if $|t - \hat \tau| < \delta$.
\end{enumerate}
\end{lemma}

By combining the above results, we have the following lemma, which is useful in our proof of Theorem \ref{Thm_Main}.

\begin{lemma}\label{Lem_ExactOne}
Suppose $p(t)$ satisfies the conditions in Theorem \ref{Thm_Main}
and $\tau_0$ satisfies the properties in Lemma \ref{Lem_PenaltyProperty2}.
Let $h(\theta, \mu)$ and $t^*(\theta, \mu)$ be as defined in Lemma \ref{Lem_MinUnique}
and Lemma \ref{Lem_MinUnique_1} respectively for the case $q > 1$ and $q = 1$.
Then we can find $\theta$ and $\mu$ such that for any $l \geq 2$, $t_1, \ldots, t_l \in \mathbb{R}$,
\begin{equation*}
\sum _{j=1}^l p(|t_j|) + \theta \cdot \left| \sum_{j=1}^l t_j  \right|^q + \mu \cdot \left| \sum_{j=1}^l t_j - \hat \tau \right|^q  \ge h(\theta, \mu).
\end{equation*}
Moreover, let
$
\bar \delta = \min \left\{ \frac{\tau_0}{3}, \frac{t^*(\theta, \mu) - \tau_0}{2}, \frac{\tau - t^*(\theta, \mu)}{2}, 1, C_1
\right\}
$
where $C_1$ is defined in Lemma \ref{Lem_PenaltyProperty2}, then for any $\delta \in (0, \bar \delta)$, we have
\begin{equation}   \label{Eq:Lem_ExactOne_2}
\sum _{j=1}^l p(|t_j|) + \theta \cdot \left| \sum_{j=1}^l t_j  \right|^q + \mu \cdot \left| \sum_{j=1}^l t_j - \hat \tau \right|^q  < h(\theta, \mu) + \delta^2
\end{equation}
holds only if $|t_i - t^*(\theta, \mu)| < 2 \delta $ for some $i$ while $|t_j| \leq \delta$ for all $j\neq i$.
\end{lemma}

\medskip
\medskip

\begin{proof}[Proof of Theorem \ref{Thm_Main}] We present a polynomial time reduction to problem
\eqref{Eq:Problem_Main} from the 3-partition problem. For any given instance of the 3-partition problem with $\bb=(b_1, \ldots,
b_{3m})$, we consider the minimization problem
$ \min_{ x  }f( x )$
in the form of \eqref{Eq:Problem_Main} with $ x =\{x_{ij}\}, 1\leq i \leq 3m,
1 \leq j \leq m$, where
\begin{eqnarray*}
 f( x ) :=
 \sum_{j=2}^m \left|\sum_{i=1}^{3m} b_i x_{ij}-\sum_{i=1}^{3m} b_i x_{i1} \right|^q
 + \sum ^{3m} _{i=1}\left| \left( \lambda \theta \right)^{\frac 1 q}   \sum _{j=1}^m  x_{ij} \right|^q\\
 + \sum ^{3m} _{i=1}\left| \left( \lambda \mu \right)^{\frac 1 q} \left( \sum _{j=1}^m  x_{ij}- \hat \tau \right)  \right|^q
 + \lambda \sum_{i=1}^{3m} \sum_{j=1}^m p(|x_{ij}|).
\end{eqnarray*}
Note that the lower bounds $\underline \theta$, $\underline \mu$, and $\hat{\underline \mu}$
only depend on the penalty function $p(\cdot)$,
we can choose $\theta \geq \underline \theta$ and $\mu \geq \underline \mu \theta$ if $q > 1$,
or $\theta = 0$ and $\mu \geq \hat{\underline \mu}$ if $q = 1$,
such that $\left( \lambda \theta \right)^{1/q}$ and $\left( \lambda \mu
\right)^{1/q}$ are both rational numbers. Since $\hat \tau$ is also rational,
all the coefficients of $ f( x )$ are of finite size
and independent of the input size of the given 3-partition instance.
Therefore, the minimization problem $ \min_{ x  }
f( x )$ has polynomial size with respect to the given 3-partition
instance.

For any $ x $, by Lemma \ref{Lem_ExactOne},
\begin{equation} \label{Eq:P_min_value}
    \begin{aligned}
    f( x ) \geq &0 + \lambda \cdot \sum_{i=1}^{3m} \Bigg\{
\sum_{j=1}^m p(|x_{ij}|) + \theta\cdot \left| \sum _{j=1} ^m
x_{ij}\right|^q \\
&+ \mu \cdot \left| \sum _{j=1} ^m x_{ij}
-\hat\tau\right|^q \Bigg\} \geq  {3m} \lambda \cdot h(\theta, \mu).
\end{aligned}
\end{equation}

Now we claim that there exists an equitable partition to the
3-partition problem if and only if the optimal value of
$f( x )$ is smaller than ${3m}\lambda\cdot h(\theta, \mu)+\epsilon$ where $\epsilon$ is specified later. On one hand,
if $S$ can be equally partitioned into $m$ subsets, then we define
\begin{equation*}
x_{ij}= \left\{
\begin{array}{ll}
t^*(\theta, \mu)& \mbox{if $b_i$ belongs to the $j$th subset};\\
0&\mbox{otherwise}.
\end{array}\right.
\end{equation*}
It can be easily verified that these $x_{ij}$'s satisfy $f(x) = {3m} \lambda \cdot h(\theta, \mu)$. Then due to \eqref{Eq:P_min_value}, we know that these $x_{ij}$'s
provide an optimal
solution to $f( x )$ with optimal value ${3m} \lambda \cdot
h(\theta, \mu)$.

On the other hand, suppose the optimal value of $f( x )$ is
${3m} \lambda \cdot h(\theta, \mu)$, and there is a polynomial-time
algorithm that solves \eqref{Eq:Problem_Main}. Then for
\begin{eqnarray*}
\delta = \min\left\{ \frac{\tau_0 }{8\sum_{i=1}^{3m} b_i}, \bar \delta \right\}
\quad\mbox{ and }\quad
\epsilon = \min \{ \lambda \delta^2, (\tau_0 / 2)^q \}
\end{eqnarray*}
where
\begin{align*}
\bar \delta = \min \bigg\{ &\frac{\tau_0}{3}, \frac{t^*(\theta, \mu) - \tau_0}{2}, \frac{\tau - t^*(\theta, \mu)}{2},\\
&\frac{p(\tau_0 / 3) + p(2 \tau_0 / 3) - p(\tau_0)}{\tau_0 / 3} , 1 \bigg\},
\end{align*}
we are able to find a near-optimal solution $x$ such that
$f(x) < {3m}\lambda \cdot h(\theta, \mu) + \epsilon$ within a
polynomial time of $\log(1/\epsilon)$ and the size of $f(x)$, which is polynomial with respect to the size of the given
3-partition instance. Now we show that we can find an equitable
partition based on this near-optimal solution. By the definition of
$\epsilon$, $f(x) < {3m}\lambda \cdot h(\theta, \mu) +
\epsilon$ implies
\begin{equation}
    \begin{aligned}
\label{Eq:Eq:Pf:Thm_3}
&\sum_{j=1}^m p(|x_{ij}|) + \theta \left| \sum _{j=1} ^m x_{ij}\right|^q + \mu \cdot \left| \sum _{j=1} ^m x_{ij} -\tau\right|^q\\
<& h(\theta, \mu) +  \delta^2
, \quad \forall i = 1, \ldots, {3m}.
\end{aligned}
\end{equation}
According to Lemma \ref{Lem_ExactOne}, for each $i = 1, \ldots, {3m}$,
\eqref{Eq:Eq:Pf:Thm_3} implies that there exists $k$ such that
$|x_{ik} - t^*(\theta, \mu)| <2 \delta$ and $|x_{ij}| <\delta$ for
any $j \neq k$. Now let
$$
y_{ij} = \left\{
\begin{array}{c l}
t^*(\theta, \mu) & \text{if } |x_{ik} - t^*(\theta, \mu)| < 2\delta  \\
0   & \text{if } |x_{ij} | < \delta
\end{array}.
\right.
$$
We define a partition by assigning $b_i$ to the $j$th subset
$S_j$ if $y_{ij} =t^*(\theta, \mu)$. Note that this partition is
well-defined since for each $i$, by the definition of $\delta$,
there exists one and only one $y_{ik} = t^*(\theta, \mu)$ while the
others equal $0$. Now we show that this is an equitable partition.

Note that for any $j = 1, \ldots, m$, the difference between the sum of the $j$-th subset and the first subset is
\begin{align*}
\left| \sum_{S_j} b_i - \sum_{S_1} b_i \right| = &
\left|\sum_{i=1}^{3m} \frac{y_{ij}}{t^*(\theta, \mu)}\cdot  b_i
-\sum_{i=1}^{3m} \frac{y_{i1}}{t^*(\theta, \mu)  }\cdot b_i \right|\\
 = & \frac{1}{t^*(\theta, \mu)} \left|\sum_{i=1}^{3m} b_i
y_{ij}-\sum_{i=1}^{3m} b_i y_{i1} \right|.
\end{align*}
By triangle inequality, we have 
\begin{align*}
 &\left| \sum_{S_j} b_i - \sum_{S_1} b_i \right|\le \frac{1}{t^*(\theta,\mu)} \left(\sum_{i=1}^{3m}  b_i\cdot
|y_{ij}-x_{ij}|\right. \\
&\qquad+\left.\sum_{i=1}^{3m} b_i\cdot |y_{i1}-x_{i1}| +
\left|\sum_{i=1}^{3m} b_i x_{ij}-\sum_{i=1}^{3m} b_i x_{i1}
\right|\right).
\end{align*}
By the definition of $y_{ij}$, we have $|y_{ij} - x_{ij}| < 2\delta$
for any $i$, $j$. for the last term, since $f(x) < {3m}\lambda
\cdot h(\theta, \mu) + \epsilon$, we know that
\begin{align*}
\left|\sum_{i=1}^nb_ix_{ij} - \sum_{i=1}^nb_ix_{i1}\right| <
\epsilon^{1/q} \le \tau_0/2.
\end{align*}
Therefore, we have
\begin{align*}
\left| \sum_{S_j} b_i - \sum_{S_1} b_i \right| <
\frac{1}{t^*(\theta,\mu)}\left(4\delta \sum_{i=1}^nb_i +
\frac{\tau_0}{2}\right) \le 1.
\end{align*}
Now since $b_i$'s are all integers, we must have $\sum_{S_j} b_i =
\sum_{S_1} b_i$, which means that the partition is equitable.
\end{proof}


\bibliography{NPHard,reference}

\begin{thebibliography}{28}
\providecommand{\natexlab}[1]{#1}
\providecommand{\url}[1]{\texttt{#1}}
\expandafter\ifx\csname urlstyle\endcsname\relax
  \providecommand{\doi}[1]{doi: #1}\else
  \providecommand{\doi}{doi: \begingroup \urlstyle{rm}\Url}\fi

\bibitem[Amaldi \& Kann(1998)Amaldi and Kann]{amaldi1998approximability}
Amaldi, E. and Kann, V.
\newblock On the approximability of minimizing nonzero variables or unsatisfied
  relations in linear systems.
\newblock \emph{Theoretical Computer Science}, 209\penalty0 (1):\penalty0
  237--260, 1998.

\bibitem[Antoniadis \& Fan(2001)Antoniadis and
  Fan]{antoniadis2001regularization}
Antoniadis, A. and Fan, J.
\newblock Regularization of wavelet approximations.
\newblock \emph{Journal of the American Statistical Association}, 96\penalty0
  (455):\penalty0 939--967, 2001.

\bibitem[Arora et~al.(1993)Arora, Babai, Stern, and Sweedy]{arora1993hardness}
Arora, S., Babai, L., Stern, J., and Sweedy, Z.
\newblock The hardness of approximate optima in lattices, codes, and systems of
  linear equations.
\newblock In \emph{Foundations of Computer Science, 1993. Proceedings., 34th
  Annual Symposium on}, pp.\  724--733. IEEE, 1993.

\bibitem[Bian \& Chen(2014)Bian and Chen]{bian2014optimality}
Bian, W. and Chen, X.
\newblock Optimality conditions and complexity for non-lipschitz constrained
  optimization problems.
\newblock \emph{Preprint}, 2014.

\bibitem[Cameron \& Trivedi(2013)Cameron and Trivedi]{cameron2013regression}
Cameron, A.~C. and Trivedi, P.~K.
\newblock \emph{Regression analysis of count data}, volume~53.
\newblock Cambridge university press, 2013.

\bibitem[Candes et~al.(2008)Candes, Wakin, and Boyd]{candes2008enhancing}
Candes, E., Wakin, M., and Boyd, S.
\newblock Enhancing sparsity by reweighted ${L}_1$ minimization.
\newblock \emph{Journal of Fourier Analysis and Applications}, 14\penalty0
  (5-6):\penalty0 877--905, 2008.

\bibitem[Chartrand(2007)]{chartrand2007exact}
Chartrand, R.
\newblock Exact reconstruction of sparse signals via nonconvex minimization.
\newblock \emph{Signal Processing Letters, IEEE}, 14\penalty0 (10):\penalty0
  707--710, 2007.

\bibitem[Chen et~al.(2014)Chen, Ge, Wang, and Ye]{ChenX_GeD_2014}
Chen, X., Ge, D., Wang, Z., and Ye, Y.
\newblock Complexity of unconstrained ${L}_2-{L}_p$ minimization.
\newblock \emph{Mathematical Programming}, 143\penalty0 (1-2):\penalty0
  371--383, 2014.

\bibitem[Chen \& Wang(2016)Chen and Wang]{ChenY_2016}
Chen, Y. and Wang, M.
\newblock Hardness of approximation for sparse optimization with ${L}_0$ norm.
\newblock \emph{Technical Report}, 2016.

\bibitem[Davis et~al.(1997)Davis, Mallat, and Avellaneda]{davis1997adaptive}
Davis, G., Mallat, S., and Avellaneda, M.
\newblock Adaptive greedy approximations.
\newblock \emph{Constructive approximation}, 13\penalty0 (1):\penalty0 57--98,
  1997.

\bibitem[Fan \& Li(2001)Fan and Li]{fan2001variable}
Fan, J. and Li, R.
\newblock Variable selection via nonconcave penalized likelihood and its oracle
  properties.
\newblock \emph{Journal of the American Statistical Association}, 96\penalty0
  (456):\penalty0 1348--1360, 2001.

\bibitem[Fan \& Lv(2010)Fan and Lv]{fan2010selective}
Fan, J. and Lv, J.
\newblock A selective overview of variable selection in high dimensional
  feature space.
\newblock \emph{Statistica Sinica}, 20\penalty0 (1):\penalty0 101--148, 2010.

\bibitem[Fan et~al.(2014)Fan, Xue, and Zou]{FanJ_XueL_2014}
Fan, J., Xue, L., and Zou, H.
\newblock Strong oracle optimality of folded concave penalized estimation.
\newblock \emph{The Annals of Statistics}, 42\penalty0 (3):\penalty0 819--849,
  2014.

\bibitem[Fan et~al.(2015)Fan, Liu, Sun, and Zhang]{fan2015tac}
Fan, J., Liu, H., Sun, Q., and Zhang, T.
\newblock T{A}{C} for sparse learning: Simultaneous control of algorithmic
  complexity and statistical error.
\newblock \emph{arXiv preprint arXiv:1507.01037}, 2015.

\bibitem[Foster et~al.(2015)Foster, Karloff, and Thaler]{foster2015variable}
Foster, D., Karloff, H., and Thaler, J.
\newblock Variable selection is hard.
\newblock In \emph{COLT}, pp.\  696--709, 2015.

\bibitem[Frank \& Friedman(1993)Frank and Friedman]{frank1993statistical}
Frank, L.~E. and Friedman, J.~H.
\newblock A statistical view of some chemometrics regression tools.
\newblock \emph{Technometrics}, 35\penalty0 (2):\penalty0 109--135, 1993.

\bibitem[Garey \& Johnson(1978)Garey and Johnson]{garey1978strong}
Garey, M.~R. and Johnson, D.~S.
\newblock ``{S}trong''{N}{P}-completeness results: Motivation, examples, and
  implications.
\newblock \emph{Journal of the ACM (JACM)}, 25\penalty0 (3):\penalty0 499--508,
  1978.

\bibitem[Huber(1964)]{huber1964robust}
Huber, P.~J.
\newblock Robust estimation of a location parameter.
\newblock \emph{The Annals of Mathematical Statistics}, 35\penalty0
  (1):\penalty0 73--101, 1964.

\bibitem[Huo \& Chen(2010)Huo and Chen]{huo2010complexity}
Huo, X. and Chen, J.
\newblock Complexity of penalized likelihood estimation.
\newblock \emph{Journal of Statistical Computation and Simulation}, 80\penalty0
  (7):\penalty0 747--759, 2010.

\bibitem[Loh \& Wainwright(2013)Loh and Wainwright]{loh2013regularized}
Loh, P.-L. and Wainwright, M.~J.
\newblock Regularized {M}-estimators with nonconvexity: Statistical and
  algorithmic theory for local optima.
\newblock In \emph{Advances in Neural Information Processing Systems}, pp.\
  476--484, 2013.

\bibitem[Lv \& Fan(2009)Lv and Fan]{LvJ_FanY_2009}
Lv, J. and Fan, Y.
\newblock A unified approach to model selection and sparse recovery using
  regularized least squares.
\newblock \emph{The Annals of Statistics}, 37\penalty0 (6A):\penalty0
  3498--3528, 2009.

\bibitem[McCullagh(1984)]{mccullagh1984generalized}
McCullagh, P.
\newblock Generalized linear models.
\newblock \emph{European Journal of Operational Research}, 16\penalty0
  (3):\penalty0 285--292, 1984.

\bibitem[Natarajan(1995)]{natarajan1995sparse}
Natarajan, B.~K.
\newblock Sparse approximate solutions to linear systems.
\newblock \emph{SIAM journal on computing}, 24\penalty0 (2):\penalty0 227--234,
  1995.

\bibitem[Wang et~al.(2014)Wang, Liu, and Zhang]{wang2014optimal}
Wang, Z., Liu, H., and Zhang, T.
\newblock Optimal computational and statistical rates of convergence for sparse
  nonconvex learning problems.
\newblock \emph{Annals of statistics}, 42\penalty0 (6):\penalty0 2164, 2014.

\bibitem[Xue et~al.(2012)Xue, Zou, Cai, et~al.]{xue2012nonconcave}
Xue, L., Zou, H., Cai, T., et~al.
\newblock Nonconcave penalized composite conditional likelihood estimation of
  sparse ising models.
\newblock \emph{The Annals of Statistics}, 40\penalty0 (3):\penalty0
  1403--1429, 2012.

\bibitem[Zhang(2010{\natexlab{a}})]{ZhangCH_2010}
Zhang, C.-H.
\newblock Nearly unbiased variable selection under minimax concave penalty.
\newblock \emph{The Annals of Statistics}, 38\penalty0 (2):\penalty0 894--942,
  2010{\natexlab{a}}.

\bibitem[Zhang(2010{\natexlab{b}})]{ZhangT_2010}
Zhang, T.
\newblock Analysis of multi-stage convex relaxation for sparse regularization.
\newblock \emph{Journal of Machine Learning Research}, 11:\penalty0 1081--1107,
  2010{\natexlab{b}}.

\bibitem[Zhang et~al.(2014)Zhang, Wainwright, and Jordan]{zhang2014lower}
Zhang, Y., Wainwright, M.~J., and Jordan, M.~I.
\newblock Lower bounds on the performance of polynomial-time algorithms for
  sparse linear regression.
\newblock In \emph{COLT}, 2014.

\end{thebibliography}
\bibliographystyle{icml2017}

\end{document}